
\input gtmacros
\input amsnames
\input amstex
\input epsf         
\input newinsert    
%
\catcode`\@=12        
%
\input gtmonout
\volumenumber{2}
\volumeyear{1999}
\volumename{Proceedings of the Kirbyfest}
\pagenumbers{177}{199}
\papernumber{10}
\received{30 August 1999}\revised{14 October 1999}
\published{18 November 1999}
\let\\\par
\def\topmatter{\relax}
\def\endtopmatter{\maketitle}
\let\gttitle\title
\def\title#1\endtitle{\gttitle{#1}}
\let\gtauthor\author
\def\author#1\endauthor{\gtauthor{#1}}
\let\gtaddress\address
\def\address#1\endaddress{\gtaddress{#1}}
\def\affil#1\endaffil{\gtaddress{#1}}
\let\gtemail\email
\def\email#1\endemail{\gtemail{#1}}
\def\subjclass#1\endsubjclass{\primaryclass{#1}}
\let\gtkeywords\keywords
\def\keywords#1\endkeywords{\gtkeywords{#1}}
\def\heading#1\endheading{{\def\S##1{\relax}\def\\{\relax\ignorespaces}
    \section{#1}}}
\def\head#1\endhead{\heading#1\endheading}

\def\subhead#1\endsubhead{\sh{#1}}
\def\subsubhead#1\endsubsubhead{\sh{#1}}
\def\specialhead#1\endspecialhead{\sh{#1}}
\def\demo#1{\rk{#1}\ignorespaces}
\def\enddemo{\ppar}
\let\remark\demo
\def\endremark{}

\def\qed{\ifmmode\quad\sq\else\hbox{}\hfill$\sq$\par\goodbreak\rm\fi}  
\def\proclaim#1{\rk{#1}\sl\ignorespaces}
\def\endproclaim{\rm\ppar}
\def\cite#1{[#1]}
\newcount\itemnumber
\def\roster{\items\itemnumber=1}
\def\endroster{\enditems}
\let\itemold\item
\def\item{\itemold{{\rm(\number\itemnumber)}}%
\global\advance\itemnumber by 1\ignorespaces}
\def\S{section~\ignorespaces}  
\def\date#1\enddate{\relax}
\def\thanks#1\endthanks{\relax}   
\def\dedicatory#1\enddedicatory{\relax}  
\let\footnote\plainfootnote

\def\Refs{\ppar{\large\bf References}\ppar\bgroup\leftskip=25pt
\frenchspacing\parskip=3pt plus2pt\small}       
\def\endRefs{\egroup}
\def\widestnumber#1#2{\relax}
\def\endrefitem{}
\def\refdef#1#2#3{\def#1{\leavevmode\unskip\endrefitem#2\def\endrefitem{#3}}}
\def\ref{\par}
\def\endref{\endrefitem\par\def\endrefitem{}}
\refdef\key{\noindent\llap\bgroup[}{]\ \ \egroup}
\refdef\no{\noindent\llap\bgroup[}{]\ \ \egroup}
\refdef\by{\bf}{\rm, }
\refdef\manyby{\bf}{\rm, }
\refdef\paper{\it}{\rm, }
\refdef\book{\it}{\rm, }
\refdef\jour{}{ }
\refdef\vol{}{ }
\refdef\yr{(}{) }
\refdef\ed{(}{, Editor) }
\refdef\publ{}{ }
\refdef\inbook{from: ``}{'', }
\refdef\pages{}{ }
\refdef\page{}{ }
\refdef\paperinfo{}{ }
\refdef\bookinfo{}{ }
\refdef\publaddr{}{ }
\refdef\moreref{}{ }
\refdef\finalinfo{}{ }
\refdef\eds{(}{, Editors)}
\refdef\bysame{\hbox to 3 em{\hrulefill}\thinspace,}{ }
\refdef\toappear{(to appear)}{ }
\refdef\issue{no.\ }{ }
\newcount\refnumber\refnumber=1
\def\refkey#1{\expandafter\xdef\csname cite#1\endcsname{\number\refnumber}%
\global\advance\refnumber by 1}
\def\cite#1{[\csname cite#1\endcsname]}
\def\Cite#1{\csname cite#1\endcsname}  
\def\key#1{\noindent\llap{[\csname cite#1\endcsname]\ \ }}

\refkey{B1}
\refkey{B2}
\refkey{B3}
\refkey{B4}
\refkey{BPZ}
\refkey{Bo}
\refkey{Boy}
\refkey{BZ1}
\refkey{BZ2}
\refkey{BZ3}
\refkey{BM}
\refkey{CJ}
\refkey{CGLS}
\refkey{E}
\refkey{EW}
\refkey{FS}
\refkey{Ga1}
\refkey{Ga2}
\refkey{Ga3}
\refkey{Ga4}
\refkey{GS}
\refkey{G1}
\refkey{G2}
\refkey{G3}
\refkey{GLi}
\refkey{GL1}
\refkey{GL2}
\refkey{GL3}
\refkey{GL4}
\refkey{GL5}
\refkey{GL6}
\refkey{GW1}
\refkey{GW2}
\refkey{GW3}
\refkey{HM}
\refkey{JS}
\refkey{J}
\refkey{K}
\refkey{Me}
\refkey{Mes}
\refkey{M}
\refkey{MM}
\refkey{N}
\refkey{O}
\refkey{Q}
\refkey{S}
\refkey{Sc1}
\refkey{Sc2}
\refkey{T1}
\refkey{T2}
\refkey{W1}
\refkey{W2}
\refkey{W3}

\let\smc\small
\NoBlackBoxes

\def\F{{\Cal F}}

\def\que{{\Bbb Q}}

\def\zed{{\Bbb Z}}
\def\Int{\operatorname{Int}}
\def\lk{\operatorname{lk}}
\def\rel{\operatorname{rel}}

\def\un#1{\smallskip\noindent\underbar{#1}}
\def\strut{\hbox{\vrule height 1.3em depth.6em width0pt}}
\def\fig#1#2#3{\midinsert \cl{\epsfysize#1truein\epsfbox{#2}}
	\vglue 6pt\centerline{\smc Figure #3}\endinsert}
\topmatter
\title Small surfaces and Dehn filling\endtitle
\author Cameron McA Gordon\endauthor
\address Department of Mathematics, The University of Texas at Austin\\
Austin, TX 78712-1082, USA\endaddress
\email gordon@math.utexas.edu\endemail

\abstract
We give a summary of known results on the maximal distances between
Dehn fillings on a hyperbolic 3--manifold that yield 3--manifolds
containing a surface of non-negative Euler characteristic that is
either essential or Heegaard.
\endabstract

\asciiabstract{We give a summary of known results on the maximal distances
between Dehn fillings on a hyperbolic 3-manifold that yield 3-manifolds
containing a surface of non-negative Euler characteristic that is either
essential or Heegaard.}

\primaryclass{57M25}
\secondaryclass{57M50}

\keywords
Dehn filling, hyperbolic 3--manifold, small surface
\endkeywords

\asciikeywords{Dehn filling, hyperbolic 3-manifold, small surface}

\endtopmatter
\cl{\small\it Dedicated to Rob Kirby on the occasion of his 60th birthday}

\document
\sectionnumber=-1
\head Introduction\endhead

By a {\it small\/} surface we mean one with non-negative Euler characteristic,
ie a sphere, disk, annulus or torus.
In this paper we give a survey of the results that are known on the distances
between Dehn fillings on a hyperbolic 3--manifold that yield 3--manifolds
containing small surfaces that are either essential or Heegaard.
We also give some new examples in this context.

In Section 1 we describe the role of small surfaces in the theory of
3--manifolds, and in Section~2 we summarize known results on the distances
$\Delta$ between Dehn fillings on a hyperbolic 3--manifold $M$ that create
such surfaces.
Section~3 discusses the question of how many manifolds $M$ realize the
various maximal values of $\Delta$, while Section~4 considers the situation
where the manifold $M$ is {\it large\/} in the sense of Wu \cite{W3}.
Finally, in Section~5 we consider the values of $\Delta$ for fillings on
a hyperbolic manifold $M$ with $k$ torus boundary components, as $k$
increases.

I would like to thank John Luecke and Alan Reid for useful conversations.
I would also like to thank the referee for his helpful comments; in particular
for pointing out a gap in the original proof of Theorem~5.1 and for
suggesting a considerable improvement to Theorem~3.4.

The author is partially supported by NSF grant DMS 9626550.

\head Small surfaces and 3--manifolds \endhead

The importance of small surfaces in the theory of 3--manifolds is well known.
For example, every 3--manifold (for convenience we shall assume that all
3--manifolds are compact and orientable) can be decomposed into canonical
pieces by cutting it up along such surfaces.

For spheres, this is due to Kneser \cite{K} (see also Milnor \cite{M}),
and goes as follows.
If a 3--manifold $M$ contains a sphere $S$ which does not bound a ball in $M$,
then $S$ is {\it essential\/} and $M$ is {\it reducible\/}.
Otherwise, $M$ is {\it irreducible\/}.
Then any oriented 3--manifold $M$ can be expressed as a connected sum
$M_1 \# \ldots \# M_n$, where each $M_i$ is either irreducible or
homeomorphic to $S^2\times S^1$.
Furthermore, if we insist that no $M_i$ is the 3--sphere, then the summands
$M_i$ are unique up to orientation-preserving homeomorphism.

Turning to disks,
a properly embedded disk $D$ in a 3--manifold $M$ is
said to be {\it essential\/} if
$\partial D$ does not bound a disk in $\partial M$.
If $M$ contains such a disk, ie, if $\partial M$ is compressible, then
$M$ is {\it boundary reducible\/}; otherwise $M$ is
{\it boundary irreducible\/}.
Then we have the following statement about essential disks in a 3--manifold,
proved by Bonahon in \cite{Bo}:
In any  irreducible 3--manifold $M$, if $W$ is a maximal (up to isotopy)
disjoint union of compression bodies on the components of $\partial M$, then
$W$ is unique up to isotopy, any essential disk in $M$ can be isotoped
$(\rel \partial)$ into $W$, and $\overline{M-W}$ is irreducible and
boundary irreducible.
Note that $\overline{M-W}$ is obtained from $M$ by cutting $M$ along a
collection  of essential disks that is maximal in the appropriate sense.

Now, let us say that a
connected, orientable, properly embedded surface $F$, not a sphere or
disk, in a 3--manifold $M$
is {\it essential\/} if it is incompressible and not parallel to a
subsurface of $\partial M$.
With this definition, an essential surface may be boundary compressible.
However, if $F$ is an essential annulus and $M$ is irreducible and boundary
irreducible, then $F$ is boundary incompressible.

Then, in an irreducible, boundary irreducible 3--manifold $M$, there
is a canonical (up to isotopy) collection $\F$ of disjoint essential annuli
and tori, such that each component of $M$ cut along $\F$ is either a
Seifert fiber space, an $I$--bundle over a surface, or a
3--manifold that contains no essential annulus or torus.
This is the JSJ--decomposition of $M$, due to Jaco and Shalen \cite{JS}
and Johannson \cite{J}.

Following Wu \cite{W3}, let us call a 3--manifold {\it simple\/} if it
contains no essential sphere, disk, annulus or torus.
Then Thurston has shown \cite{T1}, \cite{T2}
that a 3--manifold $M$ with non-empty
boundary (other than $B^3$) is simple if and only if it is {\it hyperbolic\/},
in the sense that $M$ with its boundary tori removed has a complete
hyperbolic structure with totally geodesic boundary.

For  closed 3--manifolds $M$, if $\pi_1(M)$ is finite then
Thurston's geometrization conjecture \cite{T1}, \cite{T2}
asserts that $M$ has a spherical structure.
Equivalently, $M$ is either $S^3$, a lens space, or a Seifert fiber space
of type $S^2 (p_1,p_2,p_3)$ with $\frac1{p_1} + \frac1{p_2} + \frac1{p_3} >1$.
(We shall say that a Seifert fiber space is {\it of type\/}
$F(p_1,p_2,\ldots,p_n)$ if it has base surface $F$ and $n$ singular fibers
with multiplicities $p_1,p_2,\ldots,p_n$.)
Note that $S^3$ contains a Heegaard sphere, while a lens space contains
a Heegaard torus.
For closed 3--manifolds $M$ with infinite fundamental group, there are two
cases.
If $\pi_1(M)$ has no $\zed\times\zed$ subgroup, then the geometrization
conjecture says that $M$ is hyperbolic.
If $\pi_1(M)$ does have a $\zed\times\zed$ subgroup, then by work of
Mess \cite{Mes}, Scott \cite{Sc1}, \cite{Sc2}, and, ultimately, Casson and
Jungreis \cite{CJ} and Gabai \cite{Ga4}, $M$  either contains an
essential torus or is a
Seifert fiber space of type $S^2(p_1,p_2,p_3)$.

Summarizing, we may say that if a 3--manifold is not hyperbolic then it
either
\roster
\itemold{\rm(1)} contains an essential sphere, disk, annulus or torus; or
\itemold{\rm(2)} contains a Heegaard sphere or torus; or
\itemold{\rm(3)} is a Seifert fiber space of type $S^2 (p_1,p_2,p_3)$; or
\itemold{\rm(4)} is a counterexample to the geometrization conjecture.
\endroster

\head Distances between small surface Dehn fillings\endhead

Recall that if $M$ is a 3--manifold with a torus boundary component $T_0$,
and $\alpha$ is a {\it slope\/} (the isotopy class of an essential
unoriented simple closed curve) on $T_0$, then the manifold obtained by
{\it $\alpha$--Dehn filling\/} on $M$ is $M(\alpha) = M\cup V$, where $V$
is a solid torus, glued to $M$ along $T_0$ in such a way that $\alpha$
bounds a disk in $V$.
If $M$ is hyperbolic, then
the set of {\it exceptional\/} slopes
$E(M) = \{\alpha :M(\alpha)$ is not hyperbolic$\}$
is finite \cite{T1}, \cite{T2},
and we are interested in obtaining universal upper
bounds on the size of $E(M)$.
Note that if $\alpha\in E(M)$ then $M(\alpha)$ satisfies (1), (2), (3) or
(4) above.
Here we shall focus on (1) and (2), in other words, where $M(\alpha)$
contains a small surface that is either essential or Heegaard.
(For results on case~(3), see Boyer's survey article \cite{Boy} and
references therein, and also \cite{BZ3}.)

Following Wu \cite{W3}, let us say that a 3--manifold is of {\it type\/}
$S$, $D$, $A$ or $T$ if it contains an essential sphere, disk, annulus
or torus.
Let us also say that it is of {\it type\/} $S^H$ or $T^H$ if it contains
a Heegaard sphere or torus.
Recall that the {\it distance\/} $\Delta (\alpha_1,\alpha_2)$ between two
slopes on a torus is their minimal geometric intersection number.
Then, for $X_i\in \{S,D,A,T,S^H,T^H\}$ we define

$\qquad\Delta (X_1,X_2) =
\max \{\Delta (\alpha_1,\alpha_2) : \text{there is a hyperbolic 3--manifold
$M$}$\nl
$\phantom{\qquad\Delta (X_1,X_2) =m}\text{and
slopes $\alpha_1,\alpha_2$ on a torus component of $\partial M$}$\nl
$\phantom{\qquad\Delta (X_1,X_2) =m}\text{such that $M(\alpha_i)$ is of type $X_i$, $i=1,2\}$.}$
\medskip

The numbers $\Delta (X_1,X_2)$ are now known in almost all cases,
and are summarized in Table~2.1.

\midinsert\small
\vbox{
$$\vbox{\offinterlineskip\halign{
\enspace$\strut#$\enspace
&\vrule#&\enspace$#$\enspace
&\vrule#&\hfill\enspace$#$\hfill\enspace
&\vrule#&\hfill\enspace$#$\hfill\enspace
&\vrule#&\hfill\enspace$#$\hfill\enspace
&\vrule#&\hfill\enspace$#$\hfill\enspace
&\vrule#&\hfill\enspace$#$\enspace\hfill
&\vrule#&\hfill\enspace$#$\enspace\hfill&\vrule#\cr
&&S^{\ \,}&&D^{\ \,}&&A^{\ \,}&&T^{\ \,}&&S^H&&T^H&\cr
\noalign{\hrule}
S&&1&&0&&2&&3&&?&&1&\cr
\noalign{\hrule}
D&& &&1&&2&&2&&\text{--}&&\text{--}&\cr
\noalign{\hrule}
A&& && &&5&&5&&\text{--}&&\text{--}&\cr
\noalign{\hrule}
T&& && && &&8&&2&&?&\cr
\noalign{\hrule}
S^H&& && && && &&0&&1&\cr
\noalign{\hrule}
T^H&& && && && && &&1&\cr
\noalign{\hrule}
}}$$
\centerline{\smc Table 2.1\qua $\Delta (X_1,X_2)$}
}
\endinsert

(The entries $\Delta (X_1,X_2)$ for $X_1 =D$ or $A$ and $X_2 = S^H$ or $T^H$
are blank because the first case applies only to manifolds with boundary,
while the second case applies only to closed manifolds.)

The upper bounds in the various cases indicated in Table~2.1 are due to the
following.
$(S,S)$: Gordon and Luecke \cite{GL3};
$(S,D)$: Scharlemann \cite{S};
$(S,A)$: Wu \cite{W3};
$(S,T)$: Oh \cite{O}, Qiu \cite{Q}, and Wu \cite{W3};
$(S,T^H)$: Boyer and Zhang \cite{BZ2};
$(D,D)$: Wu \cite{W1};
$(D,A)$: Gordon and Wu \cite{GW2};
$(D,T)$: Gordon and Luecke \cite{GL6};
$(A,A)$, $(A,T)$, and $(T,T)$: Gordon \cite{G1};
$(T,S^H)$: Gordon and Luecke \cite{GL4};
$(S^H,S^H)$: Gordon and Luecke \cite{GL2};
$(S^H,T^H)$ and $(T^H,T^H)$: Culler, Gordon, Luecke
and Shalen \cite{CGLS}.

References for the existence of examples realizing these upper bounds
are as follows:

\un{$(S,S)$}\qua
An example of a hyperbolic 3--manifold, with two torus boundary components,
having a pair of reducible Dehn fillings at distance~1, is given by
Gordon and Litherland in \cite{GLi}.
By doing suitable Dehn filling along the other boundary component one obtains
infinitely many hyperbolic 3--manifolds with a single torus boundary
component, having reducible fillings at distance~1.
Infinitely many such examples with two torus boundary components are given
by Eudave-Mu\~noz and Wu in \cite{EW}.

\un{$(S,A)$, $(D,A)$ and $(D,T)$}\qua
An example of a hyperbolic 3--manifold $M$, with two torus boundary
components, with Dehn fillings $M(\alpha_1)$, $M(\alpha_2)$ such that
$M(\alpha_1)$ is reducible and boundary reducible, $M(\alpha_2)$ is
annular and toroidal, and $\Delta (\alpha_1,\alpha_2)=2$, is given by
Hayashi and Motegi in [\Cite{HM}; \S12].
Infinitely many such examples are constructed by Eudave-Mu\~noz and Wu
in \cite{EW}.

\un{$(S,T)$ and $(S,T^H)$}\qua
Boyer and Zhang point out in \cite{BZ1} and [\Cite{BZ2}; Example~7.8], that
the hyperbolic 3--manifold $M=W(6)$, obtained by 6--Dehn filling (using the
usual
meridian--latitude slope co-ordinates) on the exterior $W$ of the Whitehead
link, has the property that $M(1)$ is reducible, $M(4)$ is toroidal, and
$M(\infty)$ is the lens space $L(6,1)$.
Infinitely many such hyperbolic 3--manifolds $M$ are given by Eudave-Mu\~noz
and Wu in [\Cite{EW}; Lemma~4.1 and Theorem~4.2]; ie, each $M$ has Dehn
fillings $M(\alpha_1)$, $M(\alpha_2)$, $M(\alpha_3)$ such that
$M(\alpha_1)$ is reducible, $M(\alpha_2)$ is toroidal, $M(\alpha_3)$ is a
lens space, $\Delta (\alpha_1,\alpha_2) =3$, and $\Delta(\alpha_1,\alpha_3)$
$(=\Delta (\alpha_2,\alpha_3))=1$.

\un{$(D,D)$}\qua
Infinitely many examples of hyperbolic knots in a solid torus, with a
non-trivial Dehn surgery yielding a solid torus, have been given by
Berge \cite{B1} and \cite{Ga2}.

\un{$(A,A)$ and $(A,T)$}\qua
Miyazaki and Motegi \cite{MM} and, independently, Gordon and Wu \cite{GW1},
have shown that the exterior $M$ of the Whitehead sister link has a pair
of Dehn fillings $M(\alpha_1)$, $M(\alpha_2)$, each of which is annular
and toroidal, with $\Delta (\alpha_1,\alpha_2)=5$.

\un{$(T,T)$}\qua
Thurston has shown \cite{T1} that if $M$ is the exterior of the figure
eight knot then $M(4)$ and $M(-4)$ are toroidal.

\un{$(T,S^H)$}\qua
Infinitely many examples of hyperbolic knots in $S^3$ with half-integral
toroidal Dehn surgeries are given by Eudave-Mu\~noz in \cite{E}.

\un{$(S^H,T^H)$ and $(T^H,T^H)$}\qua
Infinitely many hyperbolic knots in $S^3$ with lens space surgeries are
described by Fintushel and Stern in \cite{FS}.
A general construction of such knots is given by Berge in \cite{B2},
who has subsequently shown \cite{B3} that the knots listed in \cite{B2}
are the only ones obtainable in this way.
He has also suggested \cite{B2} that any knot in $S^3$ with a lens space
surgery might be of this form.

There is a (unique) hyperbolic knot $K$ in $S^1\times D^2$ with two
non-trivial surgeries which yield $S^1\times D^2$; see \cite{B1}.
Under an unknotted embedding of $S^1\times D^2$ in $S^3$ with $n$
meridional twists, the image of $K$ is a hyperbolic knot $K_n$ in $S^3$
with two lens space surgeries; see \cite{B1}.
(The simplest example of this kind is the $(-2,3,7)$ pretzel knot, which is
one of the knots constructed in \cite{FS}.)
Hence there are infinitely many hyperbolic 3--manifolds $M$ with Dehn
fillings $M(\alpha_1)$, $M(\alpha_2)$, $M(\alpha_3)$ such that
$M(\alpha_1) \cong S^3$, $M(\alpha_2)$ and $M(\alpha_3)$ are lens spaces,
and $\Delta (\alpha_1,\alpha_2) = \Delta(\alpha_1,\alpha_3) =
\Delta (\alpha_2 ,\alpha_3) =1$.
\medskip

We see that only two values of $\Delta (X_1,X_2)$ are unknown, namely:
$\Delta (S,S^H)$ and $\Delta (T,T^H)$.
The conjectured values are $-\infty$ and 3, and the best bounds to date
are 1 \cite{GL1} and 5 \cite{G3}, respectively.

The assertion that $\Delta (S,S^H) =-\infty$ says that no Dehn surgery
on a hyperbolic knot in $S^3$ gives a reducible manifold.
This would follow from the

\proclaim{Cabling Conjecture}
{\rm (Gonz\'alez-Acu\~na and Short \cite{GS})}\qua
If Dehn surgery on a non-trivial knot $K$ in $S^3$ gives a reducible
manifold then $K$ is a cable knot.
\endproclaim

\noindent
(Here, it is convenient to regard a torus knot as a cable of the unknot.)

In fact, the cabling conjecture and the assertion $\Delta (S,S^H)=-\infty$
are equivalent, since Scharlemann has shown \cite{S} that the former is
true for satellite knots.

Regarding $\Delta (T,T^H)$, the figure eight sister manifold $M$ has slopes
$\alpha_1,\alpha_2$ on $\partial M$ such that $M(\alpha_1)$ is toroidal,
$M(\alpha_2)$ is the lens space $L(5,1)$, and $\Delta (\alpha_1,\alpha_2)=3$
\cite{BPZ}.
In fact, there are infinitely many such hyperbolic manifolds $M$, and also
infinitely many such $M$ where $M(\alpha_2)$ is the lens space $L(7,2)$;
see Section~3.
On the other hand, it is shown in \cite{G3} that $\Delta (T,T^H)\le 5$.
Presumably $\Delta (T,T^H)=3$: there is nothing in the argument of
\cite{G3} to suggest that the bound of 5 obtained there is best possible,
while 4 is not a Fibonacci number.

\proclaim{Question 2.1}
Is there a hyperbolic manifold with a toroidal filling and a lens space
filling at distance $4$ or $5$?
\endproclaim

\head The manifolds realizing $\Delta (X_1,X_2)$\endhead

Having determined $\Delta (X_1,X_2)$, one can ask about the manifolds $M$
that have fillings realizing $\Delta (X_1,X_2)$.
Regarding the number of such manifolds,  we have

\proclaim{Theorem 3.1}
In the cases where $\Delta (X_1,X_2)$ is known, there are infinitely many
hyperbolic manifolds $M$ realizing $\Delta (X_1,X_2)$, except when
$(X_1,X_2) = (A,A),(A,T)$ or $(T,T)$.
\endproclaim

This is well known when $\Delta (X_1,X_2)=0$.
References in the other cases are given in Section~2 above.

Turning to the exceptional cases $(A,A)$, $(A,T)$ and $(T,T)$, the first
two are simultaneously described in the following theorem.
(Here, and in Theorem 3.3, $\Delta$ denotes $\Delta (\alpha_1,\alpha_2)$.)

\proclaim{Theorem 3.2}{\rm (Gordon--Wu \cite{GW1}, \cite{GW3})}\qua
Let $M$ be a hyperbolic 3--manifold such that $M(\alpha_1)$ is annular and
$M(\alpha_2)$ is annular (toroidal).
Then there are:
\roster
\itemold{\rm(1)} exactly one such manifold with $\Delta =5$;
\itemold{\rm(2)} exactly two such manifolds with $\Delta =4$; and
\itemold{\rm(3)} infinitely many such manifolds with $\Delta =3$.
\endroster
\endproclaim

The manifolds in (1) and (2)
are the same in both the annular and the toroidal case.
They are:
in (1), the exterior of the Whitehead sister (or $(-2,3,8)$ pretzel) link,
and in (2), the exteriors of the Whitehead link and the 2--bridge link
associated with the rational number 3/10.

Although the statements in Theorem~3.2 are identical in both cases $(A,A)$
and $(A,T)$, the proofs are necessarily quite different.

The next theorem describes the case $(T,T)$.

\proclaim{Theorem 3.3}{\rm (Gordon \cite{G1})}\qua
Let $M$ be a hyperbolic 3--manifold such that $M(\alpha_1)$ and $M(\alpha_2)$
are toroidal.
Then there are:
\roster
\itemold{\rm(1)} exactly two such manifolds with $\Delta =8$;
\itemold{\rm(2)} exactly one such manifold with $\Delta =7$;
\itemold{\rm(3)} exactly one such manifold with $\Delta =6$; and
\itemold{\rm(4)} infinitely many such manifolds with $\Delta =5$.
\endroster
\endproclaim

Here the manifolds in (1), (2) and (3)
are all Dehn fillings on the exterior $W$
of the Whitehead link.
Specifically, (using the usual meridian--latitude slope co-ordinates)
they are: in (1), $W(1)$ and $W(-5)$
(these are the figure eight knot exterior and the figure eight sister
manifold), in (2), $W(-5/2)$, and in (3), $W(2)$.

Of the two cases where $\Delta (X_1,X_2)$ is not known, namely
$(X_1,X_2) = (S,S^H)$ and $(T,T^H)$, recall that it is expected that there
are no examples at all realizing $(S,S^H)$.
For the other case, $(T,T^H)$, there are no examples known with $\Delta>3$.
However, the following theorem says that there are infinitely many
examples with $\Delta =3$.

\proclaim{Theorem 3.4}
For any integer $m>0$ there are
infinitely many hyperbolic 3--manifolds $M$ with Dehn fillings
$M(\alpha_1)$, $M(\alpha_2)$ such that $M(\alpha_1)$ is toroidal,
$M(\alpha_2)$ is the lens space $L(6m\pm 1,3m\mp 1)$,
and $\Delta (\alpha_1,\alpha_2)=3$.
\endproclaim

\demo{Proof}
We will construct these manifolds by suitably modifying the examples of
hyperbolic manifolds with toroidal and reducible fillings at distance 3 given
by Eudave-Mu\~noz and Wu in [\Cite{EW}; \S4].

For $p,q\in\zed$ let $T_{p,q}$ be the tangle in the 3--ball $S^3-\Int B$
shown in Figure~3.1, where $\boxed{\ n\ }$ denotes $n$ positive half-twists,
if $n\ge0$, and $|n|$ negative half-twists, if $n<0$; this is obtained
from the tangle $T_p$ shown in [\Cite{EW}; Figure~4.1(a)]
by adding $q$ horizontal half-twists beneath the $p$ vertical half-twists.
Let $T_{p,q}(r)$ be the knot or link obtained by inserting into the 3--ball
$B$ the rational tangle parametrized (in the usual way) by $r\in\que\cup
\{\infty\}$.
Let $M_{p,q}$ be the 2--fold branched covering of $T_{p,q}$.
Thus $\partial M_{p,q}$ is a torus, and $M_{p,q}(r)$ is the 2--fold
branched covering of $T_{p,q}(r)$.

\fig{1.6}{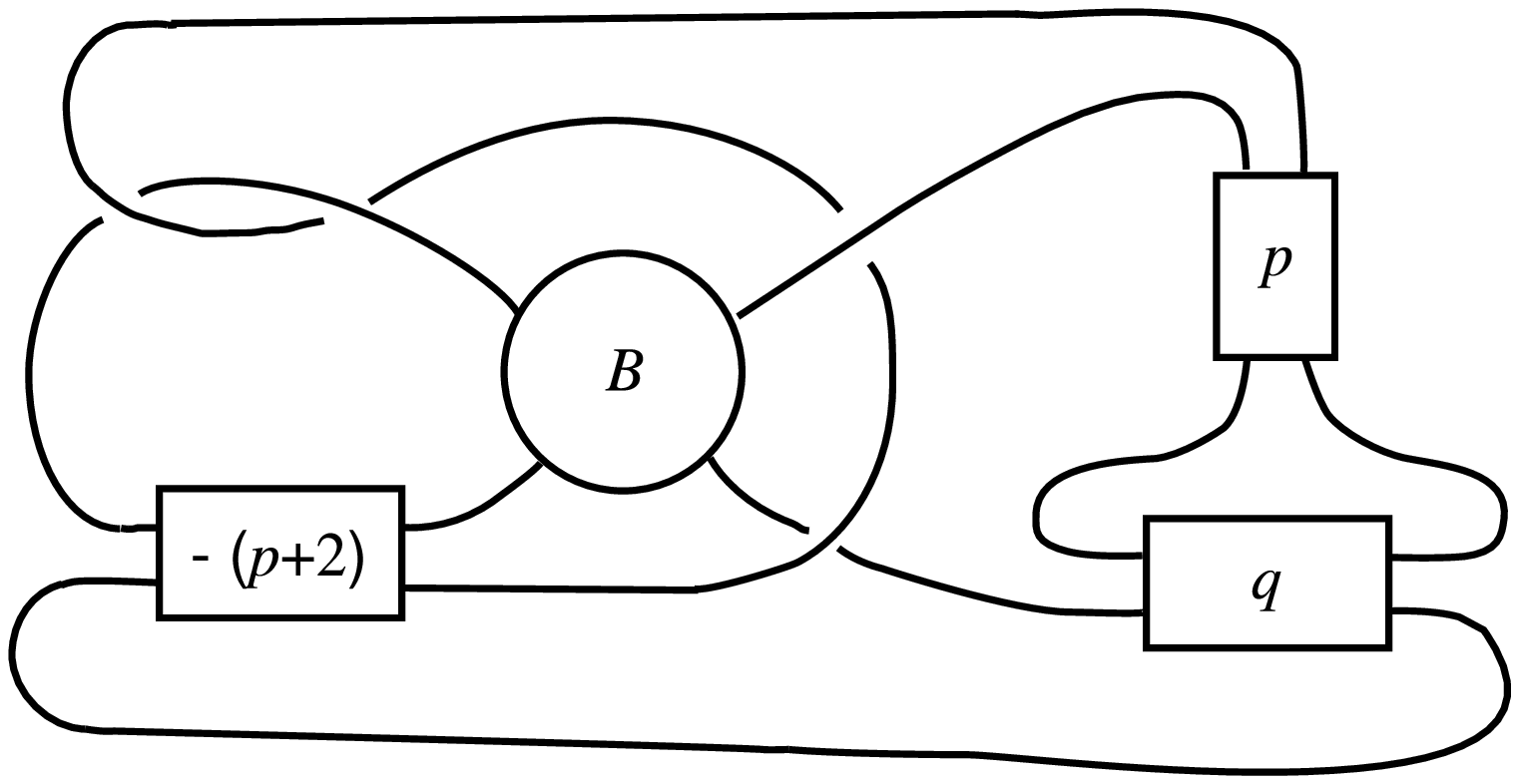}{3.1}

Assume that $p\ge3$ and $q\ne0$.
Then, as in [\Cite{EW}; Proof of Lemma~4.1], $M_{p,q}(\infty)$
is a non Seifert fibered, irreducible, toroidal manifold, $M_{p,q}(0)$ is the
2--fold branched cover of the 2--bridge knot corresponding to the rational
number $1/(-(p+3)+1/(-(p+1)+1/q))$, ie, the
lens space $L((p+3)(q(p-1)-1)+q,q(p-1)+1)$,
and $M_{p,q}(1)$ and $M_{p,q}(1/2)$ are
Seifert fiber spaces of type $S^2(p_1,p_2,p_3)$.

Also, $T_{p,q}(1/3)$ is the knot $K_q$ shown in Figure~3.2;
compare [\Cite{EW}; Figure~4.1(f)].
Thus $K_q$ is the 2--bridge knot corresponding to the rational number
$1/(-2+1/(-q+1/3)) = (1-3q)/(6q+1)$.
Hence
$M_{p,q}(1/3)$ is (up to orientation) the lens space $L(6q+1, 3q-1)$.
Setting $m=|q|$ gives the lens spaces described in the theorem.
Note that $\Delta (\infty,1/3)=3$.

\fig{1.7}{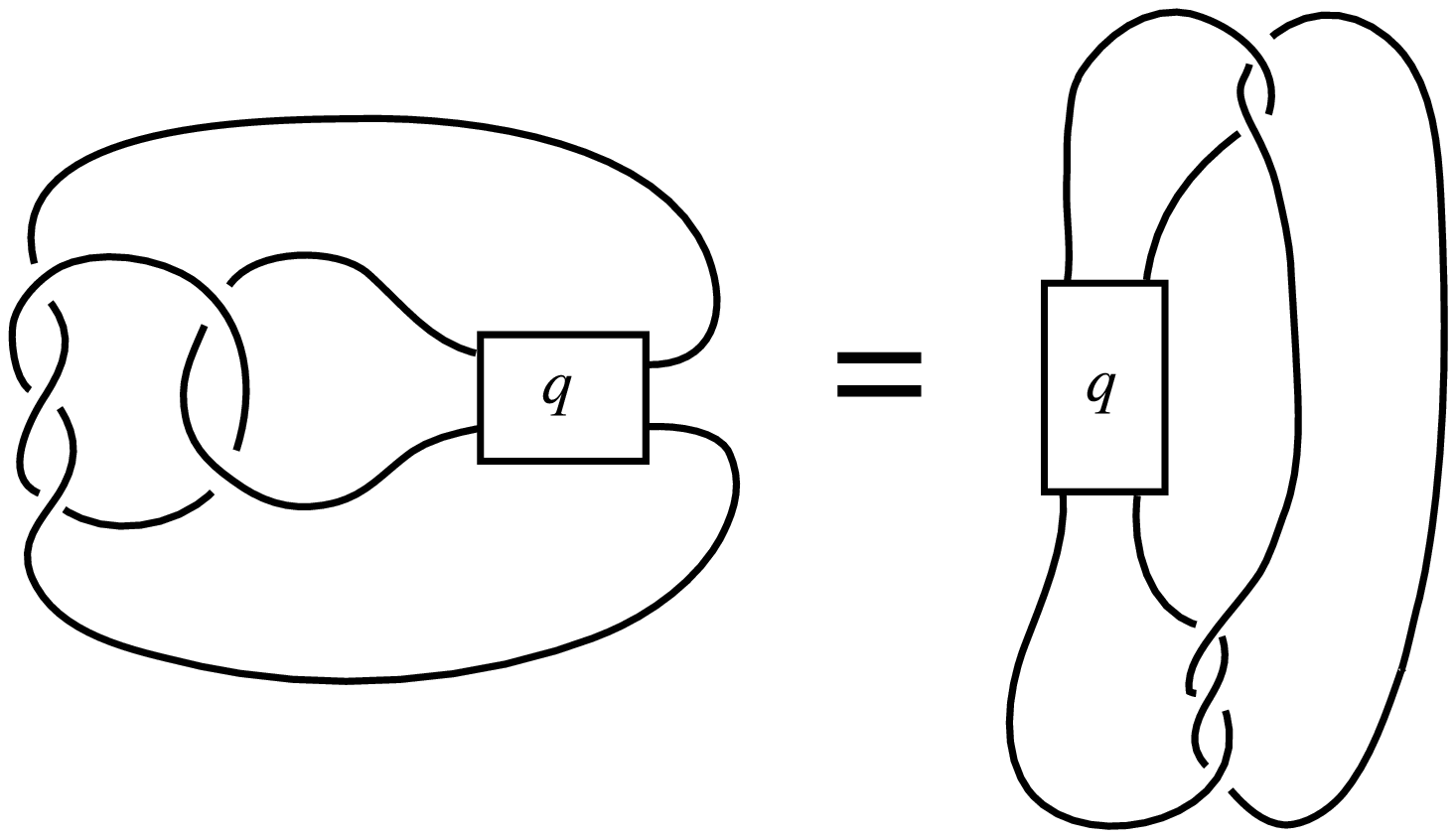}{3.2}

It remains to show that for any $q$ ($\ne0$)  there are infinitely many
distinct hyperbolic manifolds of the form $M_{p,q}$.
But the proof given by Eudave-Mu\~noz and Wu of the corresponding
assertion for their manifolds $M_p$ [\Cite{EW}; Proof of Theorem~4.2], applies
virtually unchanged in our present situation, the only modifications
necessary being to replace the reference to \cite{GL1} by one to
\cite{CGLS}, and to delete the references to \cite{GL3} and \cite{BZ2}.\qed
\enddemo

\proclaim{Question 3.1}
For which lens spaces $L$ are there infinitely many hyperbolic 3--manifolds $M$
with Dehn fillings $M(\alpha_1)$, $M(\alpha_2)$ such that $M(\alpha_1)$
is toroidal, $M(\alpha_2)$ is homeomorphic to $L$, and
$\Delta (\alpha_1,\alpha_2)=3$?
\endproclaim

\head Large Manifolds\endhead

Wu has shown \cite{W3} that for manifolds $M$ which are {\it large\/} in the
sense that $H_2 (M,\partial M-T_0) \ne0$, the bounds in Table~2.1 can often be
improved.
(Note that $M$ is {\it not\/} large if and only if
it is a $\que$--homology $S^1\times D^2$
or a $\que$--homology $T^2\times I$.)
Thus we define (for $X_i \in \{S,D,A,T\}$)

$\qua\Delta^*(X_1,X_2) = \max \{\Delta (\alpha_1,\alpha_2): \text{there is a
large hyperbolic 3--manifold $M$}$\nl
$\phantom{\qua\Delta^*(X_1,X_2) =m}\text{and slopes $\alpha_1,\alpha_2$ on a torus component of $\partial M$}$\nl
$\phantom{\qua\Delta^*(X_1,X_2) =m}\text{such that $M(\alpha_i)$ is of type $X_i$, $i=1,2\}$.}$
\medskip

(It is clear that if $M$ is large then $M(\alpha)$ can never contain
a Heegaard sphere or torus.)
Then the values of $\Delta^*(X_1,X_2)$ are as shown in Table~4.1.

\midinsert
\vbox{\small
$$\vbox{\offinterlineskip\halign{
\enspace$\strut#$\enspace
&\vrule#&\hfill\enspace$#$\hfill\enspace
&\vrule#&\hfill\enspace$#$\hfill\enspace
&\vrule#&\hfill\enspace$#$\hfill\enspace
&\vrule#&\hfill\enspace$#$\enspace\hfill
&\vrule#&\hfill\enspace$#$\enspace\hfill&\vrule#\cr
&&\ S\ &&\ D\ &&\ A\ &&\ T\ &\cr
\noalign{\hrule}
S&&0&&0&&1&&1&\cr
\noalign{\hrule}
D&& &&1&&1-2&&1&\cr
\noalign{\hrule}
A&& && &&4&&4&\cr
\noalign{\hrule}
T&& && && &&4-5&\cr
\noalign{\hrule}
}}$$
\centerline{\smc Table 4.1\qua $\Delta^* (X_1,X_2)$}}
\endinsert

The following are references for the fact that the relevant entries
in Table~4.1 are upper bounds for $\Delta^*(X_1,X_2)$.

\un{$(S,S)$}\qua
For manifolds  with boundary a union of tori this is due to Gabai
[\Cite{Ga1}; Corollary~2.4].
The general case follows from this by a trick due to John Luecke;
see [\Cite{W3}; Remark~4.2].

\un{$(S,D)$, $(D,D)$ and $(D,A)$}\qua
Here the upper bounds are the same as those for $\Delta (X_1,X_2)$ in
Table~2.1.

\un{$(S,A)$, $(S,T)$ and $(D,T)$}\qua
These are due to Wu [\Cite{W3}; Theorems~4.1 and 4.6].

\un{$(A,A)$ and $(A,T)$}\qua
By \cite{GW3}  and \cite{GW1} (see Theorem~3.2), the only hyperbolic manifold
with annular/annular or annular/toroidal fillings at distance 5 is the
Whitehead sister link exterior, which is a $\que$--homology $T^2\times I$.

\un{$(T,T)$}\qua
By \cite{G1}
(see Theorem~3.3), the only hyperbolic manifolds  with a pair of toroidal
fillings at distance greater than 5 are the fillings $W(1)$, $W(-5)$,
$W(-5/2)$ and $W(2)$ on the Whitehead link exterior $W$.
These are all $\que$--homology $S^1\times D^2$'s.
\medskip

References for the fact that the relevant entries in Table~4.1 are lower
bounds for $\Delta^* (X_1,X_2)$ are as follows.

\un{$(S,T)$ and $(D,T)$}\qua
In [\Cite{W3}; Example 4.7] Wu gives the example of the Borromean rings
exterior $M$, which has $M(\infty)$ reducible and boundary reducible and
$M(0)$ toroidal.

\un{$(S,A)$ and $(D,A)$}\qua
In [\Cite{W3}; Example 4.8] Wu constructs a hyperbolic  manifold $M$ whose
boundary consists of four tori, with slopes $\alpha_1$ and $\alpha_2$ such
that $M(\alpha_1)$ is reducible and boundary reducible, $M(\alpha_2)$
is annular, and $\Delta (\alpha_1,\alpha_2)=1$.

\un{$(D,D)$}\qua
Berge \cite{B4} and Gabai \cite{Ga3}
have given examples of simple manifolds $M$
with distinct slopes $\alpha_1$ and $\alpha_2$ such that $M(\alpha_i)$
is a handlebody of genus $g\ge2$, $i=1,2$.

\un{$(A,A)$, $(A,T)$ and $(T,T)$}\qua
It is shown in [\Cite{GW1}; Lemma 7.1] that the Whitehead link
exterior $M$ has fillings $M(\alpha_1)$, $M(\alpha_2)$, each of which is
annular and toroidal, with $\Delta (\alpha_1,\alpha_2)=4$.
Since the Whitehead link has linking number zero, $M$ is large.
\medskip

The two unknown values of $\Delta^*(X_1,X_2)$ in Table~4.1 give rise to
the following questions.

\proclaim{Question 4.1}
Is there a large hyperbolic manifold with a boundary reducible filling and
an annular filling at distance~2?
\endproclaim

\proclaim{Question 4.2}
Is there a large hyperbolic manifold with two toroidal fillings at
distance~5?
\endproclaim

\head Manifolds with boundary a union of tori\endhead

Restricting attention to hyperbolic 3--manifolds whose boundary components
are tori, we can consider what happens to the maximal distances between
exceptional fillings as the number of boundary components increases.
More precisely, we can define, for $X_i\in \{S,D,A,T\}$,

$\Delta^k(X_1,X_2)  = \max\{\Delta (\alpha_1,\alpha_2):
\text{there is a hyperbolic 3--manifold $M$ such that}$\nl
$\phantom{\Delta^k(X_1,X_2)  =m}\partial M\text{ is a disjoint union of $k$ tori, and slopes
$\alpha_1,\alpha_2$ on some}$\nl
$\phantom{\Delta^k(X_1,X_2)  =m}\text{component of $\partial M$,
such that $M(\alpha_i)$ is of type }X_i,\ i=1,2\}.$
\medskip

This is defined for $k\ge1$ if $X_1,X_2\in \{S,T\}$, and for $k\ge2$
otherwise.

Since a 3--manifold with more than two torus boundary components is large,
we have
$$\Delta^* (X_1,X_2)\ge \Delta^k (X_1,X_2)\ \text{ if }\ k\ge3.$$

If a 3--manifold whose boundary consists of $\ell$ tori contains an
essential disk, then it also contains an essential sphere, provided
$\ell\ge2$, and if it contains an essential annulus,
and is irreducible, then it also contains
an essential torus, provided $\ell\ge4$.
Hence
$$\gather
\Delta^k(S,X)  \ge \Delta^k(D,X),\ \text{if }\ k\ge 3\ ;\\
\Delta^k(T,X)  \ge \Delta^k(A,X),\ \text{if }\ k\ge 5\ \text{ and }\
\Delta^k(A,X) > \Delta^k(S,X).
\endgather$$

Now suppose $M$ is a hyperbolic 3--manifold with slopes $\alpha_1,\alpha_2$
on some torus component of $\partial M$ such that $M(\alpha_i)$ is of
type $X_i$, where $X_i = S,D,A$ or $T$, $i=1,2$.
Let $F_i$ be the corresponding essential surface in $M(\alpha_i)$, $i=1,2$.
If there is a torus component $T$ of $\partial M$ which does not meet
$F_1$ or $F_2$, then known results imply that there are infinitely many
slopes $\beta$ on $T$ such that $F_i$ remains essential in
$M(\alpha_i)(\beta)$, $i=1,2$.
Since $M(\beta)$ is hyperbolic for all but finitely many $\beta$, there are
infinitely many slopes $\beta$ such that $M(\beta)$ is hyperbolic and
$M(\beta)(\alpha_i)$ is of type $X_i$, $i=1,2$.
Thus
$$\Delta^{k-1}(X_1,X_2) \ge \Delta^k (X_1,X_2),$$
provided $k$ is large enough that there is guaranteed to be a boundary
component which misses $F_1$ and $F_2$; this depends on the pair $X_1,X_2$.

The values of $\Delta^2 (X_1,X_2)$ and $\Delta^3(X_1,X_2)$ are shown
in Tables~5.1 and 5.2.

\midinsert
\centerline{\small\raggedright\vtop{\hsize=2.5truein
\vbox{
$$\vbox{\offinterlineskip\halign{
\enspace$\strut#$\enspace
&\vrule#&\hfill\enspace$#$\hfill\enspace
&\vrule#&\hfill\enspace$#$\hfill\enspace
&\vrule#&\hfill\enspace$#$\hfill\enspace
&\vrule#&\hfill\enspace$#$\enspace\hfill
&\vrule#&\hfill\enspace$#$\enspace\hfill&\vrule#\cr
&&\ S\ &&\ D\ &&\ A\ &&\ T\ &\cr
\noalign{\hrule}
S&&1&&0&&2&&2-3&\cr
\noalign{\hrule}
D&& &&1&&2&&2&\cr
\noalign{\hrule}
A&& && &&5&&5&\cr
\noalign{\hrule}
T&& && && &&5&\cr
\noalign{\hrule}
}}$$}
\centerline{\smc Table 5.1\qua $\Delta^2 (X_1,X_2)$}}
\qquad
\vtop{\hsize=2.5truein
\vbox{
$$\vbox{\offinterlineskip\halign{
\enspace$\strut#$\enspace
&\vrule#&\hfill\enspace$#$\hfill\enspace
&\vrule#&\hfill\enspace$#$\hfill\enspace
&\vrule#&\hfill\enspace$#$\hfill\enspace
&\vrule#&\hfill\enspace$#$\hfill\enspace
&\vrule#&\hfill\enspace$#$\hfill\enspace
&\vrule#\cr
&&\ S\ &&\ D\ &&\ A\ &&\ T\ &\cr
\noalign{\hrule}
S&&0&&0&&1&&1&\cr
\noalign{\hrule}
D&& &&0&&1&&1&\cr
\noalign{\hrule}
A&& && &&3&&3&\cr
\noalign{\hrule}
T&& && && &&3-5&\cr
\noalign{\hrule}
}}$$}
\centerline{\smc Table 5.2\qua $\Delta^3 (X_1,X_2)$}}}
\endinsert

The upper bounds for $\Delta^2 (X_1,X_2)$ in Table~5.1 are the same as
the upper bounds for $\Delta(X_1,X_2)$ in Table~2.1, except for $(T,T)$.
This case follows from \cite{G1} (see Theorem~3.3), since the
manifolds listed there with a pair of toroidal fillings at distance greater
than 5 all have a single boundary component.

References for examples realizing the (lower) bounds in Table~5.1 are among
those listed for Table~2.1 in Section~2, ie,
$(S,S)$: \cite{GLi}, \cite{EW};
$(S,T)$, $(D,A)$ and $(D,T)$: \cite{HM}, \cite{EW};
$(D,D)$: \cite{B1}, \cite{Ga2};
$(A,A)$, $(A,T)$ and $(T,T)$: \cite{MM}, \cite{GW1}.

Turning to $\Delta^3 (X_1,X_2)$, the upper bounds are the same as those
for $\Delta^*(X_1,X_2)$ (see Table~4.1), except in the cases
$(D,A)$, $(A,A)$ and $(A,T)$.
For $(D,A)$, we have $\Delta^3 (D,A) \le \Delta^3 (S,A)\le1$, while the
facts that $\Delta^3(A,A)\le3$ and $\Delta^3(A,T)\le 3$
follow from \cite{GW3}
and \cite{GW1} respectively; see Theorem~3.2.

References for examples realizing the lower bounds in Table~5.2 are as
follows.

\un{$(S,T)$ and $(D,T)$}\qua
[\Cite{W3}; Example 4.7]; see Section~4 above.

\un{$(S,A)$}\qua
Let $M$ be the hyperbolic manifold constructed by Wu
in [\Cite{W3}; Example~4.8], with four torus boundary components, and slopes
$\alpha_1,\alpha_2$ (on $T_0$, say) such that $M(\alpha_1)$ is reducible,
$M(\alpha_2)$ is annular, and $\Delta(\alpha_1,\alpha_2)=1$.
By doing a suitable Dehn filling on the boundary component $T_1$ which is
neither $T_0$ nor either of the components containing the boundary
components of the annulus in $M(\alpha_2)$, we get a hyperbolic 3--manifold
$M'$ with three torus boundary components, such that $M'(\alpha_1)$ is
reducible and $M'(\alpha_2)$ is annular.
Another example is given in Theorem~5.1 below.

\un{$(D,A)$}\qua
See Theorem 5.1. (Note that although in Wu's example [\Cite{W3}; Example~4.7]
$M(\alpha_1)$ is also boundary reducible, it is $T_1$ that is compressible
in $M(\alpha_1)$, so we cannot use the argument given above in the case
$(S,A)$ to conclude that $\Delta^3 (D,A) =1$.)

\un{$(A,A)$, $(A,T)$ and $(T,T)$}\qua
In [\Cite{GW1}; Section 7] is described a hyperbolic 3--man\-if\-old $M$, called
the {\it magic manifold\/}, which is the exterior of a certain 3--component
link in $S^3$ and has Dehn fillings $M(\alpha_1)$, $M(\alpha_2)$, each of
which is annular and toroidal, with ${\Delta (\alpha_1,\alpha_2)=3}$.
\medskip

The following theorem shows that $\Delta^3 (D,A)=1$.

\proclaim{Theorem 5.1}
There exists a hyperbolic 3--component link in $S^3$ whose exterior $M$
has Dehn fillings $M(\alpha_1),M(\alpha_2)$ such that $M(\alpha_1)$ is
boundary reducible, $M(\alpha_2)$ is annular, and $\Delta (\alpha_1,
\alpha_2) = 1$.
\endproclaim

\demo{Proof}
Let $L= K_1 \cup K_2 \cup K_3$ be the 3--component link illustrated in
Figure~5.1.
Let $M$ be the exterior of $L$.
\enddemo

\midinsert{\epsfysize=2.4truein\centerline{\epsfbox{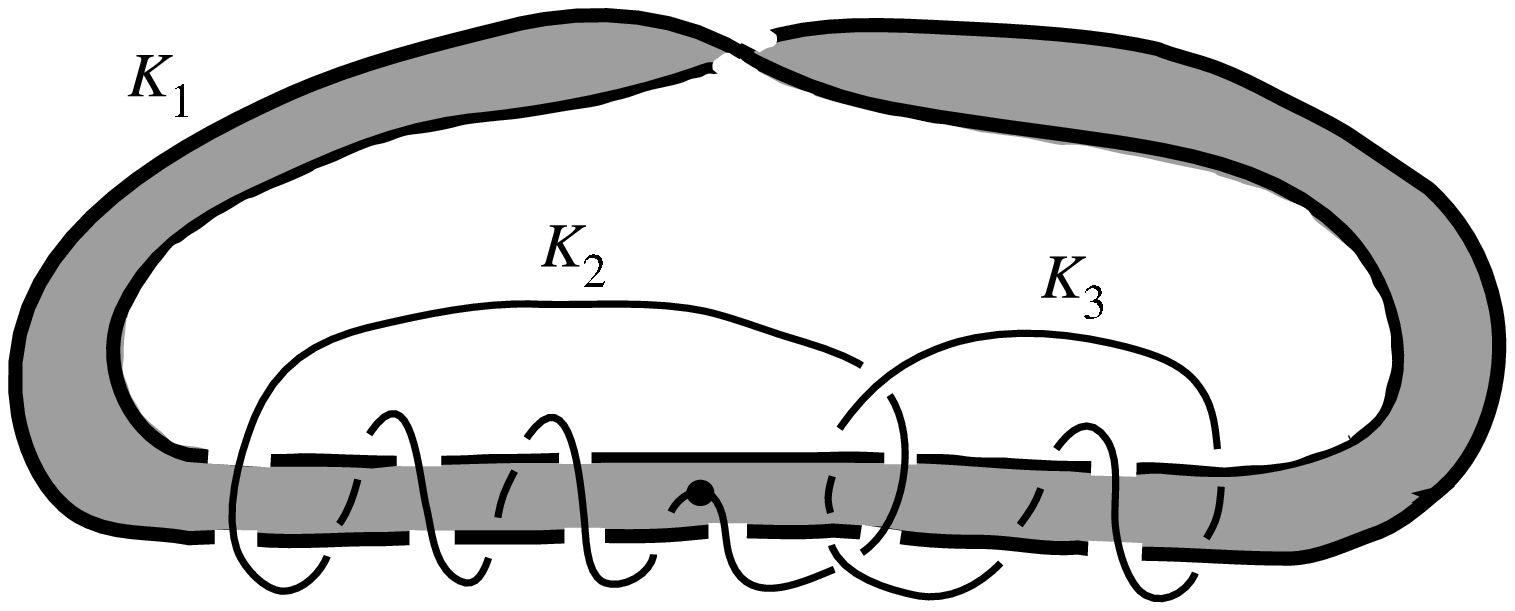}}
\vskip-1truein
\centerline{\smc Figure 5.1}}\endinsert

\proclaim{Claim}
$M$ is hyperbolic.
\endproclaim

\demo{Proof}
First, since (with appropriate orientations) we have linking numbers
$\lk (K_1,K_2)=5$, $\lk(K_1,K_3)=2$, $M$ is irreducible.
Second, it follows easily from \cite{BM}, again by considering linking
numbers, that $M$ is not a Seifert fiber space.
Hence it suffices to show that $M$ is atoroidal.

So let $T$ be an essential torus in $M$.
We see from Figure~5.1 that $K_1$ bounds a M\"obius band $B$
that is punctured once by $K_2$ and is disjoint from $K_3$.
This gives rise to a once-punctured M\"obius band $F$ in $M$.
By an isotopy of $T$, we may suppose that $T$ intersects $F$ transversely
in a finite disjoint union of simple closed curves, each being orientation
preserving and essential in $F$.
Hence we can choose an orientation reversing curve $C$ in $F$ such that
$C\cap T = \emptyset$.
Up to isotopy in $F$, there are two possibilities for $C$ (because of the
puncture), but in each case we see from Figure~5.1 that the link
$L' = C\cup K_2\cup K_3$ has a connected, prime, alternating diagram,
and is not a $(2,q)$ torus link,
and hence by \cite{Me}, is hyperbolic.
It follows that $T$ is either
\roster
\itemold{\rm(i)} compressible in $S^3- L'$; or
\itemold{\rm(ii)} parallel in $S^3-L'$ to $\partial N(C)$; or
\itemold{\rm(iii)} parallel in $S^3 - L'$ to $\partial N(K_2)$; or
\itemold{\rm(iv)} parallel in $S^3 - L'$ to $\partial N(K_3)$.
\endroster

In case (i), let $D$ be a compressing disk for $T$ in $S^3- L'$.
Then $T$ bounds a solid torus $V$ in $S^3$ containing $D$.
Since $T$ is incompressible in $S^3-L$, $D$ must meet $K_1$.
Hence $K_1\subset V$.
We now distinguish two subcases:
(a)~$K_1$ is not contained in a ball in $V$;
and
(b)~$K_1$ is contained in a ball in $V$.

In subcase (a),
since $K_1$ is unknotted in $S^3$, it follows that $V$ is also, and hence,
since $T$ is incompressible in $S^3-L$, we must have
$K_2$ or $K_3\subset S^3-V$.
If any component of $L'$ were contained in $V$, then it would lie in a
ball in $V$, and so $L'$ would be a split link.
Hence $L'\subset S^3-V$.
But $K_1\cup C$ is a Hopf link, and so $K_1$ is a core of $V$, contradicting
the essentiality of $T$ in $M$.

In subcase (b), first note  that since each of $C$, $K_2$ and $K_3$ has
non-zero linking number with $K_1$, we must have $C\cup K_2\cup K_3\subset V$,
and hence $V$ is knotted in $S^3$.
Now consider $T\cap B= T\cap F$; any component of $T\cap B$ either bounds a
disk in $B$ containing the point $K_2\cap B$, or is parallel in $B$ to $K_1$.
If there are components of the first type, let $\gamma$ be one that is
innermost in $B$; thus $\gamma$ bounds a disk $E$ in $B$ which meets $K_2$ in
a single point and has interior disjoint from $T$.
If $\gamma$ were inessential on $T$, then we would get a 2--sphere in $S^3$
meeting $K_2$ transversely in a single point, which is impossible.
Hence $E$ is a meridian disk of $V$.
But $D$ is a meridian disk of $V$ which misses $K_2$, so again we get a
contradiction.
It follows that each component of $T\cap B$ is parallel in $B$ to $K_1$.
If $T\cap B\ne\emptyset$, then the annulus in $B$ between $K_1$ and an
outermost component $\gamma$ of $T\cap B$ defines an isotopy of $K_1$,
fixing $K_3$, which takes $K_1$ to $\gamma$.
But since the meridian disk $D$ of $V$ misses $K_3$, $K_3$ lies in a ball
in $V$, and hence $\lk (\gamma,K_3)=0$.
Since $\lk (K_1,K_3)=2$, this is a contradiction.

We therefore have $T\cap B=\emptyset$.
Thus $B\subset V$, and $T$ is an essential torus in $S^3 - \text{Int }N(B\cup
K_2\cup K_3)$.
Now $B\cup K_2\cup K_3$ collapses to the graph $\Gamma\subset S^3$
shown in Figure~5.2, and $S^3-\text{Int }N(\Gamma)$ is homeomorphic to
the exterior of the tangle $t$ in $B^3$ shown in Figure~5.3.
Since $t$ is not a split tangle, $\partial B^3-t$ is incompressible in
$B^3-t$.
(To see that $t$ is not split, observe that if it were, it would be a trivial
2--string tangle together with a meridional linking circle of one of the
components.
Hence any 2--component link, with each component unknotted, obtained by
capping off $(B^3,t)$ with a trivial tangle, would be a Hopf link.
But joining the $N$ and $E$, and $S$ and $W$, arc endpoints of $t$ in the
obvious way gives the 2--bridge link corresponding to the rational number
$5/18$.)
Also, two copies of $(B^3,t)$ may be glued together along their boundaries
so as to get a link in $S^3$ that has a connected, prime, alternating
diagram.
By \cite{Me}, the exterior of this link is atoroidal, and hence the exterior
of $t$ in $B^3$ is also atoroidal.
This contradiction completes the proof of subcase~(b), and hence of case~(i).

\midinsert \cl{\hbox to 2.5truein{\vbox
{\hsize=2.5truein\epsfysize=1.5truein\centerline{\epsfbox{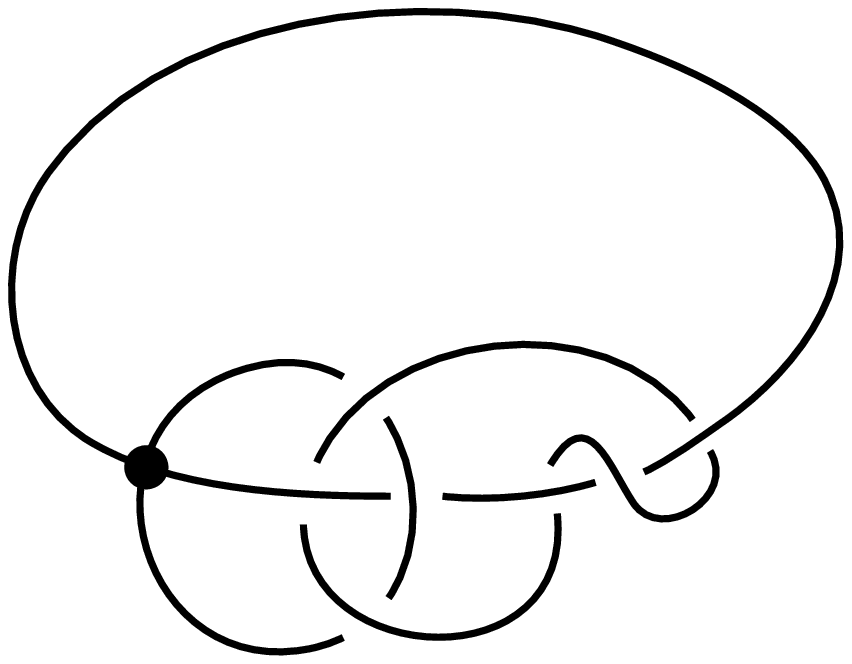}}
\vglue6pt
\centerline{\smc Figure 5.2}}}
\noindent\hbox to 2.5truein{\vbox
{\hsize=2.5truein\epsfysize=1.5truein\centerline{\epsfbox{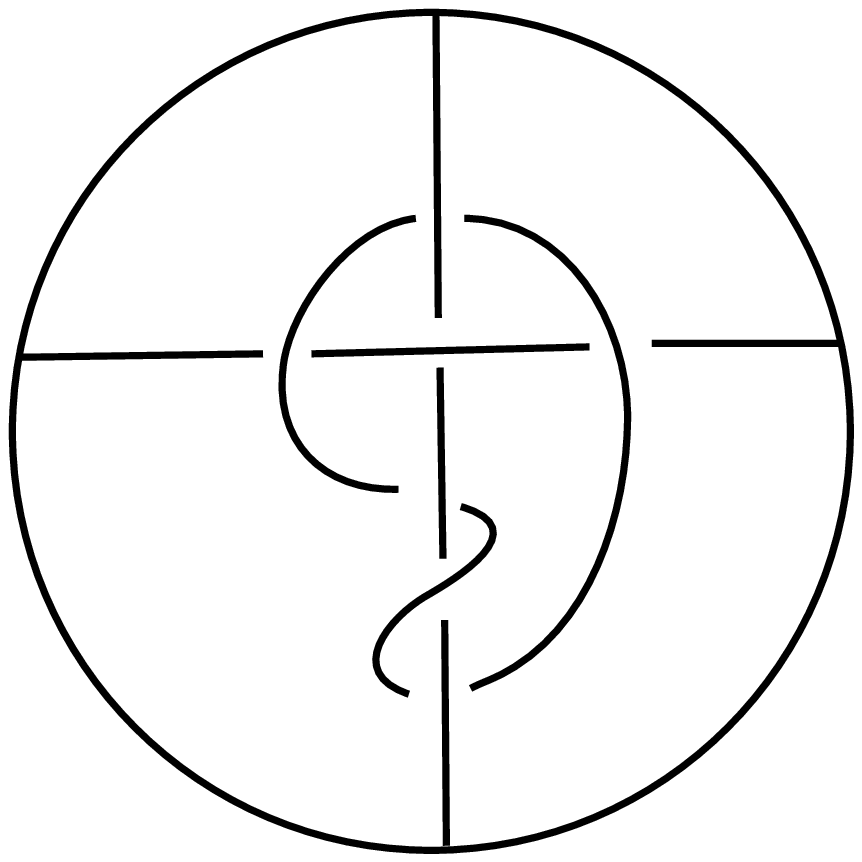}}
\vglue6pt
\centerline{\smc Figure 5.3}}}}
\endinsert

In case (ii), $T$ bounds a solid torus $V$ in $S^3$ with $C$ as a core,
and $K_1\subset V$.
Hence $\lk (K_1,K_2)=5$ is a multiple of $\lk (C,K_2)=2$ or 3, a
contradiction.
Similarly, in case (iii) we get that $\lk (K_1,C)=1$ is a multiple of
$\lk (K_2,C)=2$ or 3, and in case (iv), that $\lk (K_1,K_2)=5$ is a
multiple of $\lk (K_3,K_2)=0$.

This completes the proof of the claim.

Let $T_0$ be the boundary component of $M$ corresponding to the component
$K_1$ of $L$.
Then, since $L-K_1$ is the 2-component unlink, $M(\infty)$ is boundary
reducible.
Also, the  M\"obius band $B$ bounded by $K_1$,
which is punctured once by $K_2$, has boundary slope~2.
Hence $M(2)$ contains a M\"obius band whose boundary is a meridian of $K_2$.
Hence (see [\Cite{GLi}; Proof of Proposition~1.3, Case~(1)])
$M(2) \cong X\cup_T Q$, where $Q$ is a $(1,2)$--cable space, glued to $X$
along a torus $T$, with $Q\cap \partial M(2) =\partial N(K_2)$.
Since $M(\infty)$ is boundary reducible, $M(0)$ is irreducible, by
\cite{S}, and hence $T$ is incompressible in $X$.
Therefore $M(2)$ is annular.\qed
\enddemo

Regarding the one unknown value of $\Delta^2 (X_1,X_2)$ in Table~5.1 we
have the following question.

\proclaim{Question 5.1}
Is there a hyperbolic manifold with boundary a union of two tori, having
a reducible filling and a toroidal filling at distance~3?
\endproclaim

Similarly, the one unknown value of $\Delta^3(X_1,X_2)$ in Table~5.2 leads
to the following question.

\proclaim{Question 5.2}
Is there a hyperbolic manifold with boundary a union of three tori, having
two toroidal fillings at distance $4$ or $5$?
\endproclaim

Seeing the values in the tables for $\Delta (X_1,X_2)$, $\Delta^2 (X_1,X_2)$
and $\Delta^3 (X_1,X_2)$
decreasing leads one to ask if $\Delta^k (X_1,X_2)$
is eventually zero; equivalently, if a hyperbolic 3--manifold with $k$ torus
boundary components has at most one exceptional Dehn filling (on any given
boundary component) for $k$ sufficiently large.
However, the following two theorems show that this is not the case.

The first is essentially due to Wu \cite{W3}.

\proclaim{Theorem 5.2}{\rm (Wu \cite{W3})}\qua
For any $k\ge4$ there are infinitely many hyperbolic 3--manifolds $M$ such
that $\partial M$ consists of $k$ tori, with Dehn fillings $M(\alpha_1)$,
$M(\alpha_2)$ such that $M(\alpha_1)$ is reducible and boundary reducible,
$M(\alpha_2)$ is annular and toroidal, and $\Delta (\alpha_1,\alpha_2)=1$.
\endproclaim

\demo{Proof}
This is essentially Example 4.8 of \cite{W3}.
We simply modify Wu's construction by taking $X$ to be a simple manifold
with $\partial X$ a genus~2 surface together with $(k-4)$ tori.
Then $M= M_1\cup_PX$ is simple, with $\partial M$ consisting of $k$ tori.
Let $T_0$ be the component of $\partial M$ corresponding to $K_1$
in [\Cite{W3}; Figure~4.2].
Then $M(\infty)$ is reducible and boundary reducible, and $M(0)$ is annular.
It remains to show that $M(0)$ is toroidal.
Now $M(0)$ is irreducible (since $M$ is large and $M(\infty)$ is reducible),
and hence $M(0)$ will be toroidal
unless $k=4$ and $M(0)\cong (\text{pair of
pants}) \times S^1$.
But since $M_1(\infty)$ is reducible, $M_1(0)$ is boundary irreducible by
\cite{S}, and hence $M(0)$ contains an incompressible genus~2 surface,
so we are done.\qed
\enddemo

\remark{Remark}
Examples as in Theorem 5.2 can also be obtained by generalizing the
construction given in the proof of Theorem~5.1 to links with $k\ge 4$
components.
\endremark

It follows from Theorem 5.2 that $\Delta^k(X_1,X_2)\ge1$ for $k\ge4$, where
$X_1\in \{S,D\}$ and $X_2 \in \{A,T\}$.
The next theorem shows that for $k\ge4$, $\Delta^k(A,A)$, $\Delta^k(A,T)$
and $\Delta^k(T,T)$ are $\ge2$.

\proclaim{Theorem 5.3}
For any $k\ge4$ there exists a $k$--component hyperbolic link in
$S^2\times S^1$ whose exterior $M$ has Dehn fillings
$M(\alpha_1),M(\alpha_2)$,
each of which is annular and toroidal, with $\Delta (\alpha_1,\alpha_2)=2$.
\endproclaim

\demo{Proof}
Consider the tangle in $S^2\times I$ illustrated in Figure~5.4, consisting
of three arcs and a closed loop $K_1$
(The tangle is shown lying in the solid cylinder $D_+^2 \times I$, where
we regard $S^2$ as the union of two hemispheres $D_+^2 \cup D_-^2$.)
Gluing together the two ends $S^2\times\{0\}$ and $S^2\times \{1\}$, in
such a way that the pairs of points $\{a,a'\}$, $\{b,b'\}$ and $\{c,c'\}$
are identified, we obtain a 2--component link $L= K_1\cup K_2$ in
$S^2 \times S^1$.
For convenience we have chosen the knot $K_2$ to be the (reflection of the)
one considered by Nanyes in \cite{N}, so that we can appeal to some of
the properties of $K_2$ established there.

\fig{1.5}{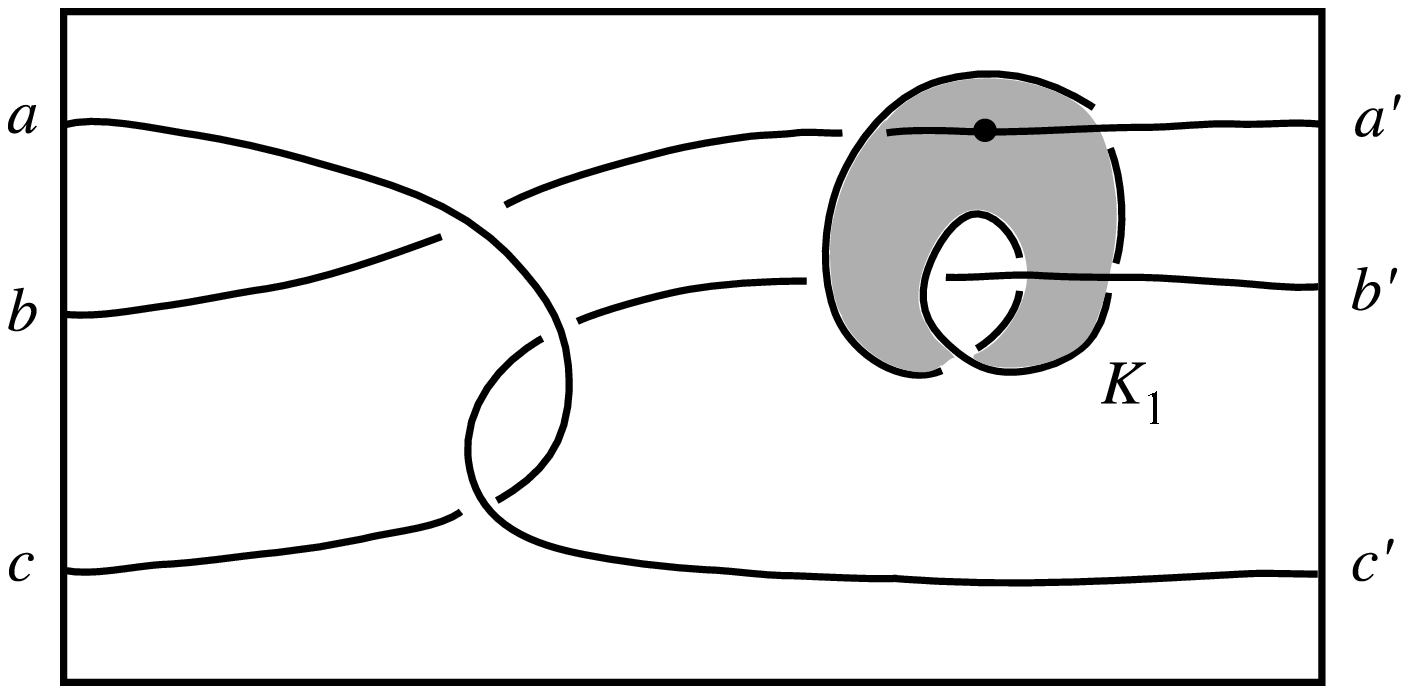}{5.4}

We see from Figure~5.4
that $K_1$ bounds a M\"obius band, with boundary slope~2,
which is punctured once by $K_2$.
Hence, doing 2--Dehn filling on the exterior of $L$ along the boundary
component $T_0$ corresponding to $K_1$, we get a manifold containing a
M\"obius band, whose boundary is a meridian of $K_2$.

Redrawing $K_1$ as in Figure 5.5, we also see that $K_1$ bounds a disk, with
boundary slope~0, which $K_2$ intersects in two points, with the same sign.
Hence 0--Dehn filling the exterior of $L$ along $T_0$ gives a manifold
that contains an annulus, whose boundary consists of two coherently
oriented meridians of $K_2$.

\fig{1.5}{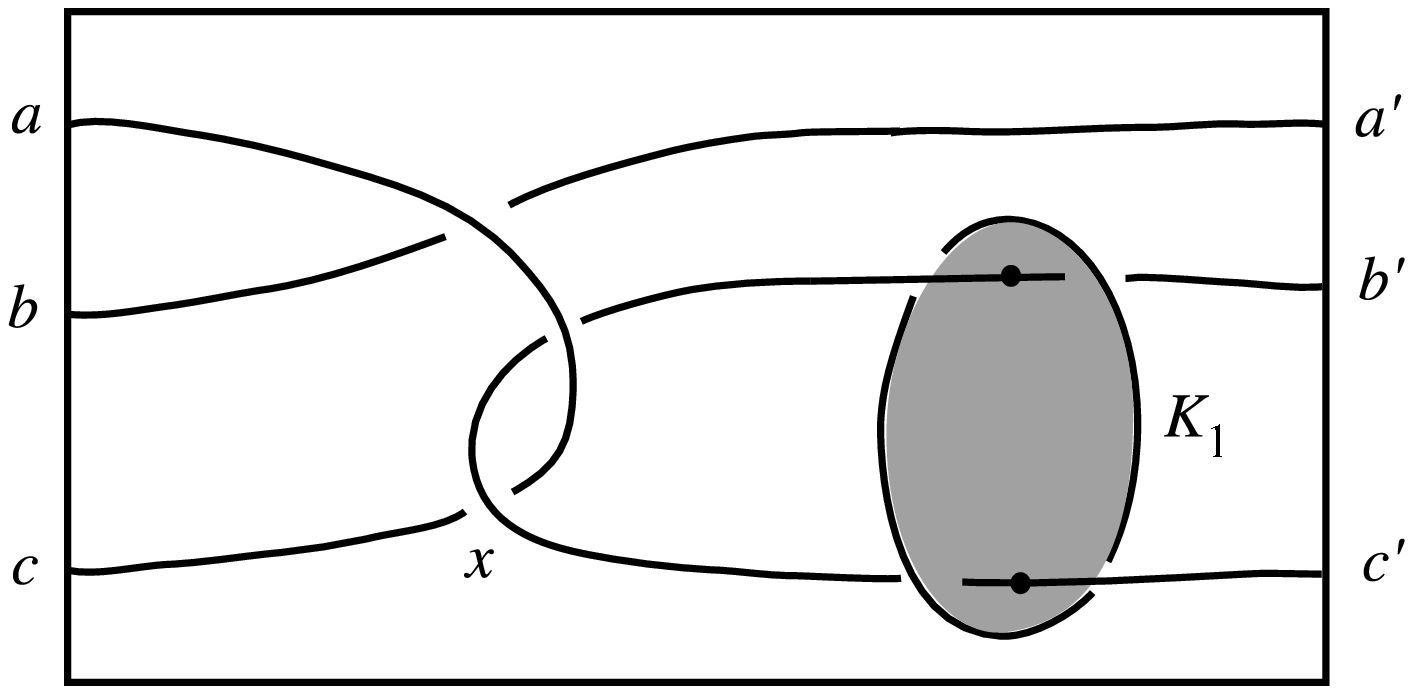}{5.5}

One can show that $L$ is hyperbolic, and the idea is to enlarge $L$ to a
$k$--component hyperbolic link $L_k$, $k\ge4$, without disturbing the
M\"obius band and annulus described above.
We do this by successively inserting $(k-2)$ additional components
$K_3,\ldots,K_k$ in a small neighborhood of the crossing $x$ indicated
in Figure~5.5, as follows.
First insert $K_3$ around $x$ as shown in Figure~5.6; then, in the same
manner,
insert $K_4$ around one of the crossings of $K_3$ with (say) $K_2$; then
insert $K_5$ around one of the crossings of $K_4$ with $K_3$ (say), etc..
Let $M$ denote the exterior of $L_k$ in $S^2\times S^1$.
Then we still have that $M(2)$ contains a M\"obius band, and $M(0)$ contains
an annulus, as described earlier.

\midinsert
\vbox{\epsfysize=1.3truein\centerline{\epsfbox{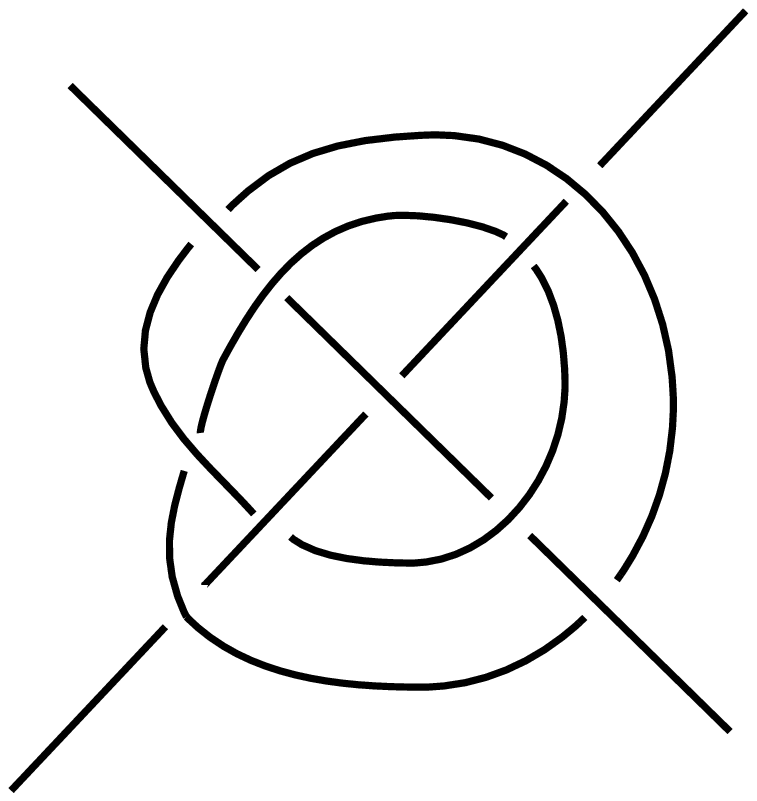}}
\vskip-.3truein
\centerline{\smc Figure 5.6}}
\endinsert

We shall show that $M$ is hyperbolic, and that $M(2)$ and $M(0)$ are
annular and toroidal.

First, let $t$ be the tangle in $S^2\times I$ that corresponds to the link
$L_k$, ie, the tangle obtained from that illustrated in Figure~5.5 by
inserting the components $K_3,\ldots,K_k$ as described above.
Let $N$ be the exterior of $t$ in $S^2\times I$.
\enddemo

\proclaim{Claim 1}
$N$ is irreducible and atoroidal.
\endproclaim

\demo{Proof}
The arc of $t$ with endpoints $a'$ and $b$ may be isotoped away from the
rest of $t$, so $N$ is homeomorphic to the exterior of the tangle $t_0$ in
$D^2\times I\cong B^3$, obtained from that shown in Figure~5.7 by
inserting $K_3,\ldots,K_k$.
Gluing two copies of $(B^3,t_0)$ along their boundaries, in such a way
that the arc endpoints $a$ and $c$ in each copy are identified with $b'$
and $c'$ respectively in the other copy, we obtain a link in $S^3$ with a
diagram that is connected, prime and alternating.
Therefore, by \cite{Me}, the exterior of this link is irreducible and
atoroidal.
Moreover, since $t_0$ is not a split tangle, $\partial B^3-t_0$ is
incompressible in $B^3-t_0$.
It follows that the exterior $N$ of $t_0$ in $B^3$ is also irreducible
and atoroidal.

This completes the proof of Claim 1.
\enddemo

\midinsert
\vbox{\epsfysize=1.8truein\centerline{\epsfbox{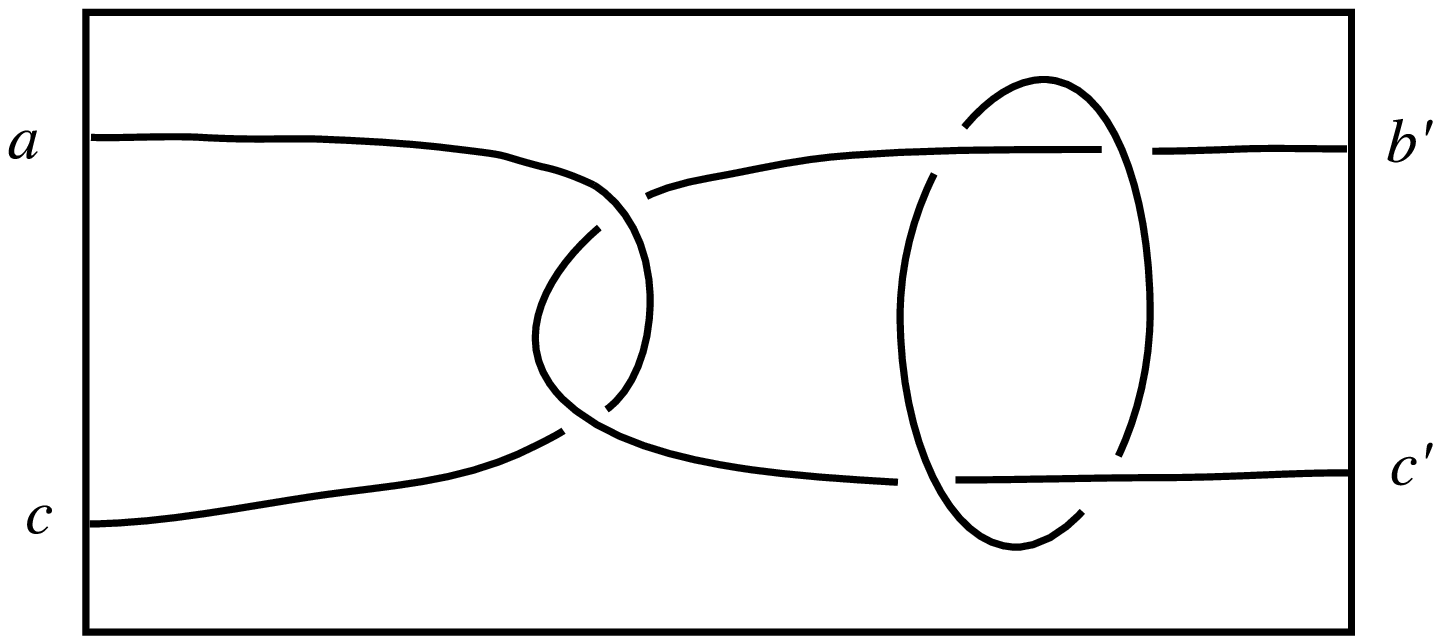}}
\vskip-.4truein
\centerline{\smc Figure 5.7}}
\endinsert

\proclaim{Claim 2}
$M$ is hyperbolic.
\endproclaim

\demo{Proof}
The two thrice-punctured spheres $P_i = S^2\times\{i\}-\text{Int }N(t)$,
$i=0,1$, are incompressible in the exterior of $t$ in $S^2\times I$, as they
are incompressible in the exterior in $S^2\times I$ of the three arcs that
make up $K_2$; see \cite{N}.
Let $P$ be the thrice-punctured sphere $P_1=P_2$ in $M$.
Since $N$ is irreducible by Claim~1, and $P_1$ and $P_2$ are
incompressible in $N$, it follows that $M$ is irreducible.

If $M$ were a Seifert fiber space, then the incompressible surface $P$
would be horizontal, which is impossible since $M$ has at least four
boundary components.

Hence it suffices to show that $M$ is atoroidal.
So let $T$ be an essential torus in $M$, which we isotop to minimize the
number of components of $T\cap P$.
Then no component of $T\cap P$ is inessential in $P$, and hence either
$T\cap P=\emptyset$,  or some component $\gamma$ of $T\cap P$ bounds a
disk $D$ in the 2--sphere $S= S^2\times \{0\} = S^2\times \{1\}$, such
that $D$ meets $K_2$ transversely in a single point and has interior
disjoint from $T$.
Now $\gamma$ is essential on $T$, otherwise we get a 2--sphere in
$S^2\times S^1$ meeting $K_2$ in a single point, contradicting \cite{N}.
Hence compressing $T$ along $D$ gives a 2--sphere $\Sigma$ meeting
$K_2$ in two points.
Since $K_2$ is locally unknotted (see \cite{N}), $\Sigma$ bounds a 3--ball
$B$ in $S^2\times S^1$ such that $(B,B\cap K_2)\cong (B^3,B^1)$.
Let $T'$ be the boundary of the solid torus $V= B-\text{Int }N(K_2)$.
Note that $T$ is obtained from $\Sigma$ by adding a tube.
If this tube lies in $B$, then $T$ is isotopic to $T'$; if it lies
outside $B$, then $T$ is parallel, in the exterior of $K_2$, to
$\partial N(K_2)$.
Since $T$ is essential in $M$, in both cases we must have $(L-K_2)\cap B
\ne \emptyset$.
If $T'$ were compressible in $V-L$, then $L$ would be a split link,
contradicting the fact that $M$ is irreducible.
Hence $T'$ is an essential torus in $M$.
Note also that $T'$ may be isotoped off $P$.
This gives an essential torus in $N$, contradicting Claim~1.

This completes the proof of Claim~2.
\enddemo

\proclaim{Claim 3}
$M(2)$ and $M(0)$ are annular and toroidal.
\endproclaim

\demo{Proof}
As observed above, $M(2)$ contains a M\"obius band.
Hence, as in the proof of Theorem~5.1,
$M(2) \cong X\cup_T Q$, where $Q$ is a $(1,2)$--cable space, glued to $X$
along a torus $T$.
If $T$ is incompressible in $X$, then $M(2)$ is annular and toroidal.
On the other hand, if $T$ compresses in $X$, then $M(2)$ is reducible.

Now consider $M(0)$.
First note that, since $M(2)$ is either annular or reducible, and
$\Delta^k (S,A) =1$, $\Delta^k (S,S)=0$, for $k\ge4$, $M(0)$ is irreducible.
Now, as we saw earlier, $M(0)$ contains an annulus $A$, whose boundary
components are coherently oriented on $\partial N(K_2)$.
It follows that $A$ is not boundary parallel in $M(0)$.
If $A$ were compressible in $M(0)$, then $M(0)$ would be boundary
reducible, and hence reducible, a contradiction.
We conclude that $M(0)$ is annular.
Now, since $\Delta^k(S,A) =1$, $k\ge4$, $M(2)$ cannot be reducible,
and hence it is annular and toroidal.

Finally, tubing $A$ along $\partial N(K_2)$ gives a Klein bottle $F$ in $M(0)$.
The boundary of a regular neighborhood of $F$ is a torus $T$ which is
essential since $M(0)$ is irreducible.
Hence $M(0)$ is toroidal.

This completes the proof of Claim~3 and hence of Theorem~5.3.\qed
\enddemo


Theorems 5.2 and 5.3 (together with Theorems~3.2 and 3.3)
show that the values of $\Delta^k(X_1,X_2)$, $k\ge4$,
are as indicated in Table~5.3.

\midinsert
\vbox{
$$\vbox{\offinterlineskip\halign{
\enspace$\strut#$\enspace
&\vrule#&\hfill\enspace$#$\hfill\enspace
&\vrule#&\hfill\enspace$#$\hfill\enspace
&\vrule#&\hfill\enspace$#$\hfill\enspace
&\vrule#&\hfill\enspace$#$\hfill\enspace
&\vrule#&\hfill\enspace$#$\hfill\enspace
&\vrule#\cr
&&\ S\ &&\ D\ &&\ A\ &&\ T\ &\cr
\noalign{\hrule}
S&&0&&0&&1&&1&\cr
\noalign{\hrule}
D&& &&0&&1&&1&\cr
\noalign{\hrule}
A&& && &&2-3&&2-3&\cr
\noalign{\hrule}
T&& && && &&2-5&\cr
\noalign{\hrule}
}}$$
\centerline{\smc Table 5.3\qua $\Delta^k (X_1,X_2)$, $k\ge4$}
}\endinsert

\proclaim{Question 5.3}
What are the values of $\Delta^k(A,A)$, $\Delta^k(A,T)$ and
$\Delta^k(T,T)$ for $k\ge4$?
\endproclaim

\Addresses\recd

\enddocument